\documentclass[11pt]{amsart}
\usepackage{amsmath,amsthm,amssymb, amscd, amsfonts}
\usepackage[all]{xy}
\usepackage{leftidx}
\usepackage[latin1]{inputenc}
\usepackage{enumerate}
\usepackage[dvipsnames]{xcolor}

\theoremstyle{plain}

\newtheorem{theorem}{Theorem}[section]
\newtheorem{corollary}[theorem]{Corollary}

\newtheorem{proposition}[theorem]{Proposition}
\newtheorem{lemma}[theorem]{Lemma}

\theoremstyle{definition}
\newtheorem{definition}[theorem]{Definition}

\theoremstyle{remark}
\newtheorem{remark}[theorem]{Remark}

\numberwithin{equation}{section}\theoremstyle{plain}

\newcommand{\I}{\mathcal{I}}

\renewcommand{\1}{\textbf{1}}

\newcommand{\A}{{\mathcal A}}
\newcommand{\B}{{\mathcal B}}
\newcommand{\C}{{\mathcal C}}
\newcommand{\D}{{\mathcal D}}
\newcommand{\F}{{\mathcal F}}

\newcommand{\Z}{{\mathcal Z}}

\newcommand{\E}{{\mathcal E}}
\newcommand{\U}{{\mathcal U}}

\newcommand{\Rep}{\operatorname{Rep}}

\newcommand\Irr{\operatorname{Irr}}
\newcommand\FPdim{\operatorname{FPdim}}

\newcommand\vect{\operatorname{Vect}}
\newcommand\svect{\operatorname{sVect}}
\newcommand\id{\operatorname{id}}
\newcommand\di{\operatorname{d}}
\newcommand\Tr{\operatorname{Tr}}

\newcommand\End{\operatorname{End}}

\newcommand\Hom{\operatorname{Hom}}
\newcommand\rk{\operatorname{rk}}

\begin{document}
\title[A fusion subcategory with maximal rank]{Fusion categories containing a fusion subcategory with maximal rank}
\author{Jingcheng Dong}
\address[Jingcheng Dong]{College of Mathematics and Statistics, Nanjing University of Information Science and Technology, Nanjing 210044, China}
\email{dongjc@njau.edu.cn}

\author{Gang Chen}
\address[Gang Chen]{School of Mathematics and Statistics, Central China Normal University, Wuhan 430079, China}
\email{chengangmath@mail.ccnu.edu.cn}

\author{Zhihua Wang}
\address[Zhihua Wang]{Department of Mathematics, Taizhou University, Taizhou 225300, China}
\email{mailzhihua@126.com}




\date{}

\begin{abstract}
In this paper, we study fusion categories which contain a proper fusion subcategory with maximal rank. They can be viewed as generalizations of near-group fusion categories. We first prove that they admit spherical structure. We then classify those which are non-degenerate or symmetric. Finally, we classify such fusion categories of rank $4$.
\end{abstract}

\maketitle

\section{Introduction}
A fusion category is a $\mathbb{C}$-linear semisimple rigid tensor category with
finitely many simple objects and finite dimensional spaces of morphisms, such that
the unit object $\1$ is simple. Classification of fusion categories is an interesting but difficult question, even some cases are  not tractable at present because it at least involves the classification of finite groups. Thus it is natural to attempt a classification of fusion categories which are ``simple" in certain sense. One class of ``simple" fusion categories is the ones admitting simple fusion rules. Two classical examples of such fusion categories are Tambara-Yamagami fusion categories \cite{Tambara1998692} and near-group fusion categories \cite{siehler2003near}.

\medbreak
In this paper, we introduce the notion of an MR fusion category, where MR stands for ``maximal rank". By definition, an MR fusion category $\C$ is a fusion category which contains a fusion subcategory $\D$ with rank $\rk(\D)=\rk(\C)-1$. Assume that $\Irr(\C) = \{\1 = X_1 ,X_2 , \cdots,X_n \}$ and $\Irr(\D) =\{X_1 ,X_2 , \cdots,X_{n - 1} \}$. It will be shown in Section \ref{sec3} that
 $$X_n \otimes X_n^\ast = \mathop \oplus
\limits_{i = 1}^{n - 1} \FPdim (X_i )X_i \oplus \kappa X_n,$$
where $\kappa$ is a non-negative integer. We denote an MR fusion category $\C$ by $\C(\D,\kappa)$ since the fusion rules of $\C$ is totally determined by the fusion subcategory $\D$ and the non-negative integer $\kappa$. It is obvious that if $\D$ is pointed then $\C$ is a near-group fusion category.  So MR fusion categories can be viewed as generalizations of near-group fusion categories.

\medbreak
In Section \ref{sec3}, we prove that the maximal rank fusion subcategory $\D$ is an integral fusion category, and hence it is tensor equivalent to the representation category of a semisimple quasi-Hopf algebra, see Proposition \ref{Prop3.2}. This implies that at moment it is impossible to classify all MR fusion categories since it involves the classification of quasi-Hopf algebras.

\medbreak
A pivotal structure on a fusion category $\C$ is an isomorphism of tensor functors $\id\to**$. In such a category, one can define the categorical dimension of an object, see Section \ref{pre}. It is conjectured in \cite[Conjecture 2.8]{etingof2005fusion} that any fusion category admits a pivotal structure.  Our  result (Theorem \ref{thm45}) shows that any MR fusion category admits a pivotal (in fact, spherical) structure.

\medbreak
In Section \ref{sec5}, we study braided MR fusion categories. Our result (Theorem \ref{integral3}) shows that a braided MR fusion category is either weakly integral, an equivariantization of a Fibonacci category $\F$, or an equivariantization of the category $\svect\boxtimes\F$. We then focus on two extreme cases: symmetric categories and non-degenerate fusion categories. If an MR fusion category is symmetric then we prove that our situation is equivalent to that in \cite{1983Characters} and hence such MR fusion categories can be completely classified, see Theorem \ref{finitegroup}. We then prove that a non-degenerate  MR fusion category is either pointed, a Fibonacci category, or a minimal extension of a slightly degenerate fusion category, see Theorem \ref{non-degenerate}.

\medbreak
The classification of fusion categories of small rank dates back to Ostrik's work \cite{ostrik2003fusion}. In that paper Ostrik classified all fusion categories of rank $2$. Up to now, although some fusion categories with additional structures have been classified, the classification of all fusion categories of a given rank stops at rank $4$. In the final part of this paper, we study MR fusion categories with rank $\leqslant 4$. Our result (Theorem \ref{rank4}) shows that all these MR fusion categories can be classified in terms of some known fusion categories.

\medbreak
In this paper, the basic theory and notions of fusion categories are given in reference \cite{etingof2005fusion,drinfeld2010braided}. The fusion categories and algebras are defined on complex number field $\mathbb{C}$.

\section{Preliminaries}\label{pre}
\subsection{Frobenius-Perron dimensions}
Let $\C$ be a fusion category and $K(\C)$ be the Grothendieck ring of $\C$.
Then the set $\Irr(\C)$ of isomorphic classes of all simple objects in $\C$ is a basis of $K(\C)$. By \cite[Theorem 8.6] {etingof2005fusion},
there is a ring homomorphism $\FPdim: K(\C) \to \mathbb{R}^+$. For any $X\in
K(\C) $, $\FPdim(X)$ is called the Frobenius-Perron dimension of $X$.
The Frobenius-Perron dimension of the fusion category $\C$ is defined as
$$\FPdim(\C)=\sum_{X \in \Irr(\C)}
(\FPdim(X))^2.$$

If the Frobenius-Perron dimension of each simple object in $\C $ is an integer then $\C$ is said to be integral. If $\FPdim(\C)$ is an integer then $\C$ is said to be weakly integral. If every simple object in the category $\C$ has dimension $1$, then $\C$ is called a pointed fusion category.

\subsection{Fusion coefficients}
Let $X$ be a simple object in $\C $ and $Y$ be any object in $\C$. We define the multiplicity of $X$ in $Y$ as $m (X, Y)=
\dim\Hom_{\C} (X, Y)$. Then $m (X, Y) = m (X^*, Y^*) $.
In addition, if $X, Y, Z\in \Irr(\C) $, then
$$m (X, Y\otimes Z)=m(Y,X \otimes Z^*) = m(Y^*,Z \otimes X^*).$$

In the literature, the numbers $m(X, Y\otimes Z)$ are often called the fusion coefficients of $\C$ and denoted by $N_{YZ}^X$.
The above property of multiplicity were first given by Nichols and Richmond in the case of semisimple Hopf algebras \cite{Nichols1996297}, and then it was extended to the fusion category settings, see \cite{dong2012frobenius} for example.

\subsection{Adjoint functors}
Let $\C $ be a fusion category. The Drinfeld center $\Z (\C) $ of $\C $ is also a fusion category. Its objects are pairs $(X, c_{-,X} ) $, where $X$ is an object in $\C$, $c_ {-,X}$ is a family of natural isomorphism $c_{V,X} :V \otimes X \overset{\sim}{\to} X \otimes V$, $\forall V \in \C$. For details on Drinfeld center, see\cite[Definition VIII4.1]{kassel1995quantum}.

Let $\F: \Z (\C) \to \C $ be the forgetful functor (i.e. $\F [(X, c_{-,X} )] = X)$ and  $\I:\C \to
\Z (\C)$ be its right adjoint functor.  Then $\I (\1) $ is a commutative algebra in $\Z (\C) $ and the unit object $1$ is a simple object in $\I (\1)$, see \cite[Lemma 3.2]{etingof2011weakly}. By \cite[Proposition 5.4]{etingof2005fusion}, we have

\[
\F(\I(V)) = \oplus _{Y \in \Irr(\C)} Y \otimes V \otimes Y^*.
\]

\subsection{Universal grading}
Let $G$ be a finite group and $\C$ be a fusion category. If $\C $ has a direct sum decomposition of Abelian subcategories

\[
\C = \mathop \oplus \limits_{g \in G} \C_g ,
\]
which  satisfies $\C_g \otimes \C_h \subseteq \C_{gh} $ and $(\C_g )^* \subseteq \C_{g ^{-1}}$, then we say that $\C$ has a $G$-grading.  We say that $\C$ is a $G$-extension of $\D$ if the grading is faithful and the trivial component is $\D$. Let $\C=\oplus_{g\in G}\C_g$ be a $G$-extension of $\D$. Then $\FPdim(\C_g)=\FPdim(\C_h)$ for all $g,h\in G$ and $\FPdim(\C)=|G|\FPdim(\D)$, see \cite[Proposition 8.20]{etingof2005fusion}.

Let $\C$ be a fusion category. The fusion subcategory $\C_{ad}$ of $\C$ generated by simple objects in $X\otimes X^{*}$ for all $X\in\Irr(\C)$ is called the adjoint subcategory of $\C$. By \cite[Corollary 3.7]{gelaki2008nilpotent} every fusion category has a canonical faithful grading $\C=\oplus_{g\in \U(\C)}\C_g$ with trivial component $\C_{ad}$. This grading is called the universal grading of $\C$ and $\U(\C)$ is called the universal grading group of $\C$.

\begin{lemma}\label{fibe_functor}
Let $\C=\oplus_{g\in G}\C_g$ be an extension of a fusion category $\D$. Assume that there exists $g\in G$ such that the rank of $\C_g$ is $1$. Then the trivial component $\D$ is equivalent to the category of finite dimensional representations of a semisimple Hopf algebra.
\end{lemma}
\begin{proof}
The tensor product of $\C$ makes $\C_g$ into a rank one module category on $\D$. By  \cite[Proposition 2.2]{ostrik2003module}, we have
a monoidal functor $F:\D\to \End(\C_g)$, where $\End(\C_g)$ is the monoidal category of endofunctors of $\C_g$. Since $\C_g$ only contains one irreducible object, $\End(\C_g)=\vect$ is the trivial category and hence  $F$ is a fiber functor on $\D$. By the reconstruction theorem for finite-dimensional Hopf algebras (see e. g. \cite[Chapter 5.3]{egno2015}), $\D$ is equivalent as a tensor category to the category of finite dimensional representations of a semisimple Hopf algebra.
\end{proof}

\subsection{Spherical fusion categories}
Let $\C$ be a fusion category and $X\in \Irr(\C)$. For a morphism $\gamma : X\to X^{**}$, we denote the trace of $\gamma$ by $\Tr_X(\gamma)$.  For any  $X\in \Irr(\C)$,  the squared norm of $X$ is defined to be the product $|X|^2= \Tr_X(\gamma)\Tr_{X^*}((\gamma^{-1})^*)$, where $\gamma : X\to X^{**}$ is an isomorphism. The global dimension of a fusion category $\C$ is the sum of squared norms of its simple objects. It is denoted by $\dim(\C)$.

\medbreak
A pivotal structure on a fusion category $\C$ is an isomorphism of tensor functors $i:\id\to**$. A pivotal fusion category is a fusion category endowed with a pivotal structure. In such a fusion category, we can define the categorical dimension of an object $X$ by $\dim(X ) = \Tr_X (i)$. Moreover, $|X|^2= |\dim(X)|^2$ in a pivotal fusion category.

A pivotal fusion category is spherical if and only if $\dim(X) = \dim(X^*)$ for all $X\in \Irr(\C)$. In a spherical category, $|X|^2= \dim(X)^2$. In particular, if $k = \mathbb{C}$, then $\dim(X)$ is (totally) real. We refer readers to \cite[Section 2.1, 2.2] {etingof2005fusion} for the facts above.

\medbreak
The pivotalization $\tilde{\C}$ of a fusion category $\C$ is given in \cite[Remark 3.1] {etingof2005fusion}.
The  simple objects of $\tilde{\C}$ are pairs $(V,\alpha)$, where $V$ is a simple object of $\C$ and $\alpha: V\simeq V^{**}$ satisfies $\alpha^{**}\alpha =\gamma$, where $\gamma:\id\to ****$ is an isomorphism of tensor functor \cite[Theorem 2.6] {etingof2005fusion}. The category $\tilde{\C}$ has a canonical pivotal structure $i:\id\to **$. In fact, $\tilde{\C}$ is spherical \cite[Proposition 5.14] {etingof2005fusion}.

For each simple object $X$ of $\C$, we have two choices of such $\alpha$. Fix one and set $(X, \alpha) = X^+$ and $(X, -\alpha) = X^-$. If we set $d=\dim(X^+)=\Tr_{X^+}(i)=\Tr_{X}(\alpha)$ then $\dim(X^-)=-d$.

\medbreak
There is an obvious tensor functor $F:\tilde{\C}\to \C$,  $F((X,\alpha))=X$. The category $\C$ is spherical if and only if the functor $F$ has a tensor section $\C\to \tilde{\C}$. Equivalently, $\tilde{\C}$ should contain a fusion subcategory such that the restriction of $F$ to this subcategory is an equivalence.

The functor $F$  maps simple objects to simple objects; that is, for any $X\in \Irr(\C)$, there are precisely two objects $X^+, X^-$ such that $F(X^{\pm}) =X$ (the choice of $X^+$ and $X^-$ is arbitrary except for $\1^+ = \1$). Moreover,
If $X\in \Irr(\C)$ is self-dual, then both $X^+$ and $X^-$ are self-dual, see \cite[Section 5.1]{2013ostrikpivotal}.

\subsection{Equivariantizations}\label{subsec24}
Let $\underline{G}$ denote the tensor category whose objects are elements of $G$, morphisms are the identity morphisms and whose tensor product is given by the multiplication in $G$. Let ${\rm\underline{Aut}}_{\otimes}\mathcal{\C}$ denote the monoidal category whose objects are tensor autoequivalences of $\C$, morphisms are isomorphisms of tensor functors and tensor product is given by the composition of functors.

\medbreak
An action of $G$ on $\mathcal{\C}$ is a tensor functor
$$T:\underline{G}\to {\rm\underline{Aut}}_{\otimes}\mathcal{\C},\quad g\mapsto T_g$$
with the isomorphism $f^X_{g,h}: T_g(X)\otimes T_h(X)\cong T_{gh}(X)$, for every $X$ in $\mathcal{\C}$.

\medbreak
Let $\mathcal{\C}$ be a fusion category with an action of $G$. Then the fusion category $\mathcal{\C}^G$, called the $G$-equivariantization of $\mathcal{\C}$, is defined as follows \cite{bruguieres2000categories,drinfeld2010braided,muger2004galois}:

(1)\quad A simple object in $\mathcal{\C}^G$ is a pair $(X, (u^X_g)_{g\in G})$, where $X\in \Irr(\C)$ is a representative of the orbits of the action of $G$ on $\Irr(\C)$, and $u^X_g: T_g(X)\to X$
is an isomorphism such that,
$$u^X_gT_g(u^X_h)= u^X_{gh}f^X_{g,h},\quad \mbox{for all\quad} g,h\in G.$$

(2)\quad A morphism $\phi: (Y,u_g^Y)\to (X,u_g^X)$ in $\mathcal{\C}^G$ is a morphism $\phi: Y\to X$ in $\mathcal{\C}$ such that $\phi u_g^Y=u_g^X\phi$, for all $g\in G$.

(3)\quad The tensor product in $\mathcal{\C}^G$ is  defined as $(Y,u_g^Y)\otimes (X,u_g^X)=(Y\otimes X, (u_g^Y\otimes u_g^X)j_g|_{Y,X})$, where $j_g|_{Y,X}:T_g(Y\otimes X)\to T_g(Y)\otimes T_g(X)$ is the isomorphism giving the monoidal structure on $T_g$.

By \cite[Proposition 4.26]{drinfeld2010braided}, we have $\FPdim(\C^G)=|G|\FPdim(\C)$.

\subsection{Braided fusion categories}
A braided fusion category $\C$ is a fusion category  with a braiding
$c_{X,Y}:X\otimes Y \overset{\sim}{\to}Y\otimes X, ~\forall X,Y\in\C$. Two objects $X,Y\in\C$ are said to centralize each other if
$c_{Y,X}c_{X,Y}=\id_{X\otimes Y}$.

Let $\D\subseteq\C$ be a fusion subcategory. Then the centralizer $\D'$ of $\D$ in $\C$ is the fusion subcategory generated by objects of $\C$ that centralize every object  of $\D$. The centralizer $\Z_2(\C):=\C'$ is called the M\"{u}ger center of $\C$.

A braided fusion category $\C$ is called non-degenerate if its M\"{u}ger center  $\Z_2(\C)=\vect$ is the trivial category. A braided fusion category $\C$ is called slightly degenerate if $\Z_2(\C)=\svect$ is the category of super vector spaces.

A braided fusion category $\C$ is symmetric if $\C=\Z_2(\C)$. For any symmetric fusion category $\C$, there exists a finite group $G$ and a central element $u\in G$ such that $\C\cong\Rep(G,u)$ \cite{deligne1990categories}, where $\Rep(G,u)$ is the category of finite-dimensional representation of $G$ and $u$ acts as parity automorphism for any $X\in\Irr(\C)$. A symmetric fusion category $\C$ is a Tannakian fusion category if $\C\cong \Rep(G)$, where the braiding is given by reflection of vector spaces.

\medbreak

A twist in a braided fusion category $\mathcal{C}$ is a  natural isomorphism $\theta:\id_\mathcal{C} \to \id_\mathcal{C}$ such that
	\begin{align}
	\theta_{X\otimes Y}=(\theta_X \otimes \theta_Y) \circ c_{Y,X} \circ c_{X,Y},
	\end{align}
for all $X,Y \in \mathcal{C}$. A twist is called a ribbon structure if $(\theta_X)^*=\theta_{X^*}$ for all $X \in \mathcal{C}$.	
A premodular fusion category  is a braided fusion category endowed with a compatible ribbon structure.

In a premodular fusion category, we can define the notion of trace for an endomorphism $\xi\in \End_{\mathcal C}(X)$, which we will denote by $\Tr(\xi)$, see \cite[Definition 4.7.1]{egno2015}. Let $\mathcal C$ be a premodular fusion category with braiding $c$. The $\mathcal{S}$-matrix $S$ of $\mathcal C$ is defined by $S:= \left(s_{X,Y}\right)_{X,Y \in \Irr(\mathcal{C})}$, where $s_{X,Y}=\Tr(c_{Y,X}c_{X,Y})$.

In a premodular fusion category $\mathcal C$, we can obtain the entries of the $S$-matrix in terms of the twists, fusion rules, and categorical dimensions via the well-known balancing equation
	\begin{align}
	s_{X,Y}=\theta_X^{-1}\theta_Y^{-1} \sum\limits_{Z \in \mathcal{O(\mathcal{C})}} N_{XY}^Z \theta_Z \di_Z,
	\end{align}
for all $X,Y \in \Irr(\mathcal{C}))$ \cite[Proposition 8.13.7]{egno2015}.	

A premodular tensor category $\mathcal C$ is said to be modular if the $S$-matrix $S$ is non-degenerate.

\section{MR fusion categories}\label{sec3}

Recall that the rank of $\C$ is the cardinality of the set $\Irr(\C)$.

\begin{definition}\label{lemma3.1}
Let $\C $ be a fusion category with rank of $n$ and $\D$ be a fusion subcategory of $\C$.
If the rank of $\D$ is $n-1$ then $\D $ is called a fusion subcategory with maximal rank.
\end{definition}

It can be seen that the fusion subcategory generated by all invertible simple objects in a near-group fusion category is a fusion subcategory with maximal rank.

\begin{proposition}\label{Prop3.2}
Let $\C$ be a fusion category with rank of $n$.
If $\D$ is a fusion subcategory with maximal rank then $\D$ is an integral fusion category. In particular,
$\D$ is the representation category of a semisimple quasi-Hopf algebra.
\end{proposition}

\begin{proof}
Let $\Irr(\C) = \{\1 = X_1 ,X_2 , \cdots,X_n \}$, $\Irr(\D) =\{X_1 ,X_2 , \cdots,X_{n - 1} \}$. Obviously, $X_n $ is self-dual, i.e., $X_n =
X_n^*$. Set $X_n \otimes X_n^* = \sum\limits_{i = 1}^n {a_i X_i } $,
where $a_i $ is the multiplicity of $X_i$ in $X_n \otimes X_n^*$

Because $\D$ is a fusion subcategory, $X_n$ can not appear in the decomposition of tensor product of any two simple objects of $D$.
Therefore, for any $1\leq i,j \le n - 1$, $m(X_n ,X_j \otimes X_i^\ast ) = m(X_j ,X_n \otimes X_i ) =0 $. From the arbitrariness of $j$, the decomposition of $X_n \otimes X_i$  only contains $X_n$ as its summands. On the other hand, the equations below
\[
m(X_i ,X_n \otimes X_n^\ast ) = m(X_i^\ast ,X_n \otimes X_n^\ast )= m(X_n
,X_i^\ast \otimes X_n ) = m(X_n ,X_n \otimes X_i ) = a_i
\]
indicates that  $X_n \otimes
X_i = a_i X_n $ is the direct sum decomposition of $X_n \otimes X_i$. By using ring homomorphism $\FPdim$ on both sides, we can get
\[
\FPdim (X_n \otimes X_i ) = \FPdim (X_n )\FPdim (X_i ) = a_i \FPdim (X_n ).
\]

This shows that $\FPdim (X_i) = a_i $ is an integer. Therefore, $\D$ is an integral fusion category. It follows
from \cite[Theorem 8.3]{etingof2005fusion} that $\D$ is the representation category of a semisimple quasi-Hopf algebra.
\end{proof}

\begin{remark}
(1)\, Let $G(\C)$ be the group generated by all $1$-dimensional simple objects of $\C$. Then the fusion rules of $\C$ shows that $g\otimes X_n=X_n$ for all $g\in G(\C)$. Hence $\C$ can not be slightly degenerate, otherwise it contradicts \cite[Proposition2.6]{etingof2011weakly}.

\medbreak
(2)\, If $a_n = 0$, then $\C = \C_0 \oplus \C_1 $ has a $\mathbb{Z}_2$-grading, of which $\C_0 =\D$, $\C_1 $ contains only one simple object $X_n$. In this case, $\D$ is the representation category of a semisimple Hopf algebra, see Lemma \ref{fibe_functor}.

\medbreak

(3)\, The above proof shows that $X_n \otimes X_n^\ast = \mathop \oplus
\limits_{i = 1}^{n - 1} \FPdim (X_i )X_i \oplus a_n X_n $. If $\D $ is a pointed fusion category, then $\C $ is obviously a near-group category. Hence fusion categories containing a fusion subcategory with maximal rank are generalizations of near-group fusion categories.
\end{remark}

\begin{definition}
A fusion category $\C$ is called an MR fusion category if $\C$ contains a fusion subcategory $\D$ with maximal rank. An MR fusion category is denoted by $\C(\D,\kappa)$, where $\kappa=a_n$.
\end{definition}

\section{Any MR fusion category is spherical}\label{sec4}
In this section, $\C=\C(\D,\kappa)$ is an MR fusion category. Assume $\Irr(\C) = \{\1 = X_1 ,X_2 , \cdots, X_n \}$ and $\Irr(\D) =\{X_1 ,X_2 , \cdots,X_{n - 1} \}$. By Proposition \ref{Prop3.2}, $\D$ is integral and hence it is pseudo-unitary and admits a unique pivotal (spherical) structure. Thus, the categorical dimensions of all simple objects are positive, and coincide with their Frobenius- Perron dimensions, see \cite[Proposition 8.23, 8.24] {etingof2005fusion}.

\medbreak
 Assume that $\tilde{\C}$ is the pivotalization of $\C$ and $F:\tilde{\C}\to\C$ is the forgetful functor. The restriction $F|_{\tilde{\D}}:\to\D$ has a tensor section $G:\D\to \tilde{\D}$ such that $G(X)=(X,\alpha)$, where $\alpha$ is deduced from the canonical pivotal structure of $\D$. Hence there exists a fusion subcategory $\tilde{\D}$ of $\tilde{\C}$ which is tensor equivalent to $\D$. It follows that we can choose $\alpha_i: V_i\simeq V_i^{**}$ such that $\{\1 = X_1^+ ,X_2^+ , \cdots, X_{n-1}^+ \}$ generates a fusion subcategory equivalent to $\D$. In this case, $d_i=\dim(X_i^+ )=\dim(X_i )=\FPdim(X_i)$ for all $1\leqslant i\leqslant n-1$.

\begin{lemma}\label{lemma41}
For all $1\leqslant i\leqslant n-1$, $X_i^+\otimes X_n^+=d_iX_n^+$, $X_i^-\otimes X_n^+=d_iX_n^-$.
\end{lemma}
\begin{proof}
Applying the forgetful functor $F$, we have $F(X_i^+\otimes X_n^+)=F(X_i^+)\otimes F(X_n^+)=X_i\otimes X_n=d_iX_n$, where the last equality follows from the fusion rules of $\C$. Hence we may assume that $X_i^+\otimes X_n^+=a_iX_n^++b_iX_n^-$, where $a_i+b_i=d_i$. On the other hand, $d_id_n=\dim(X_i^+\otimes X_n^+)=\dim(a_iX_n^++b_iX_n^-)=(a_i-b_i)d_n$. Thus  $a_i-b_i=d_i$.  It follows that $a_i=d_i,b_i=0$ and $X_i^+\otimes X_n^+=d_iX_n^+$.

\medbreak
Similarly we have  $X_i^-\otimes X_n^+=d_iX_n^-$.
\end{proof}

By Lemma \ref{lemma41}, for all $1\leqslant i\leqslant n-1$, $X_i^+$ appears in the decomposition of $X_n^+\otimes (X_n^+)^*$ with multiplicity $d_i$, however $X_i^-$ does not appear in the decomposition of $X_n^+\otimes (X_n^+)^*$. Hence we may assume that
\begin{equation}\label{decomxn}
\begin{split}
X_n^+\otimes (X_n^+)^*=\sum_{i=1}^{n-1}d_iX_i^++sX_n^++tX_n^-.
\end{split}
\end{equation}
Applying the forgetful functor $F$, we have

$$F(X_n^+\otimes (X_n^+)^*)=X_n\otimes (X_n)^*=\sum_{i=1}^{n-1}d_iX_i+(s+t)X_n.$$
This means that $s+t=\kappa$. After renaming $X_n^+$, we may assume that $s-t\geqslant 0$.

\begin{lemma}\label{lemma42}
$\dim(X_n^+)=\frac{(s-t)\pm\sqrt{(s-t)^2+4a}}{2}$, $\dim(\tilde{\C})=4a+(s-t)^2\pm(s-t)\sqrt{(s-t)^2+4a}$, where $a=\sum_{i=1}^{n-1}d_i^2=\dim(\D)$.
\end{lemma}
\begin{proof}
Considering the dimension on both sides of equation (\ref{decomxn}), we have
$$d_n^2=\sum_{i=1}^{n-1}d_i^2+(s-t)d_n=a+(s-t)d_n,$$
hence $\dim(X_n^+)=\frac{(s-t)\pm\sqrt{(s-t)^2+4a}}{2}$. It follows that
$$\dim(\tilde{\C})=2\dim(\D)+2d_n^2=4a+(s-t)^2\pm(s-t)\sqrt{(s-t)^2+4a}.$$
\end{proof}

\begin{lemma}\label{lemma43}
The fusion subcategory generated by $\{X_i^+|1\leqslant i\leqslant n\}$ is tensor equivalent to the fusion category $\C$ if one of the following holds:

 (1)\, $s=\kappa$;

 (2)\, $\tilde{\C}$ is pseudo-unitary;

 (3)\, $\sqrt{\kappa^2+4a}$ is an integer;

 (4)\, $\sqrt{(s-t)^2+4a}$  is an integer.
 \end{lemma}
\begin{proof}
(1)\, If $s=\kappa$ then $t=0$ and hence $\{X_i^+|1\leqslant i\leqslant n\}$ is tensor closed. The fusion rules show that the fusion subcategory generated by $\{X_i^+|1\leqslant i\leqslant n\}$ is tensor equivalent to the fusion category $\C$, through the forgetful functor $F$.

(2)\, If $\tilde{\C}$ is pseudo-unitary then $\dim(\tilde{\C})=\FPdim(\tilde{\C})=2\FPdim(\C)$. This means that $4a+(s-t)^2+(s-t)\sqrt{(s-t)^2+4a}=4a+\kappa^2+\kappa\sqrt{\kappa^2+4a}$. On the  other hand, $0\leqslant s-t\leqslant \kappa$. So we have $s=\kappa$. The result then follows from (1).

(3)\, If $\sqrt{\kappa^2+4a}$ is an integer then $\frac{\kappa+\sqrt{\kappa^2+4a}}{2}$ is a rational algebraic integer. As we know that a rational algebraic integer must be an integer. Therefore $\FPdim(X_n)=\frac{\kappa+\sqrt{\kappa^2+4a}}{2}$ is an integer and hence $\C$ is pseudo-unitary by \cite[Proposition 8.24] {etingof2005fusion}, and so is $\tilde{\C}$. The result then follows from (2).

(4)\,  If $\sqrt{(s-t)^2+4a}$  is an integer then $\dim(X_n^+)=\frac{(s-t)\pm\sqrt{(s-t)^2+4a}}{2}$ is an integer. Then $\tilde{\C}$ is pseudo-unitary by \cite[Lemma A.1]{HONG20101000}. The result then follows from (2).
\end{proof}

The following lemma is well-known in the theory of  algebraic integers.
\begin{lemma}\label{lemma44}
Assume that $a,b,c,d$ are integers  and $\sqrt{b},\sqrt{d}$ are not integers. Then $\frac{a+\sqrt{b}}{c+\sqrt{d}}$ is an algebraic integer if and only if $\frac{a-\sqrt{b}}{c-\sqrt{d}}$ is an algebraic integer.
 \end{lemma}

\begin{theorem}\label{thm45}
Any MR fusion category $\C$ is spherical.
 \end{theorem}
\begin{proof}
Our proof idea is to prove that $\C$ is tensor equivalent to the fusion subcategory generated by $\{X_i^+|1\leqslant i\leqslant n\}$. Then $\C$ is tensor equivalent to a full tensor subcategory of a spherical category and therefore it is spherical. By Lemma \ref{lemma43}, it suffices to prove that $s=\kappa$.

\medbreak
By \cite[Proposition 8.22] {etingof2005fusion},
$$x:=\frac{\dim(\tilde{\C})}{\FPdim(\tilde{\C})}=\frac{4a+(s-t)^2\pm(s-t)\sqrt{(s-t)^2+4a}}{4a+\kappa^2+\kappa\sqrt{\kappa^2+4a}}$$
is  an algebraic integer. By Lemma \ref{lemma44},
$$y:=\frac{4a+(s-t)^2\mp(s-t)\sqrt{(s-t)^2+4a}}{4a+\kappa^2-\kappa\sqrt{\kappa^2+4a}}$$
is also an algebraic integer. Hence
$$xy=\frac{16a^2+4a(s-t)^2}{16+4a\kappa^2}$$
is an algebraic integer. It follows that it is an integer since it is a rational number. On the other, $0\leq s-t\leq \kappa$ implies that $s-t= \kappa$. Together with the fact $s+t=\kappa$, we get $s=\kappa$. This completes the proof.
\end{proof}

\section{Braided MR fusion categories}\label{sec5}
\subsection{General results}
Recall that a Fibonacci category is a rank $2$  modular category of Frobenius-Perron dimension $\frac{5+\sqrt{5}}{2}$. It is known that Fibonacci categories fall into $2$ equivalence classes and both of them can be realized using the quantum group $U_q(sl_2)$ for $q=\sqrt[10]{1}$, see \cite{ostrik2003fusion}.
\begin{lemma}\label{centralizer}
Let $\C=\C(\D,\kappa)$ be a braided MR fusion category. Assume that $\C$ is not symmetric. Then $\D=\C'$ if and only if $\D'=\C$.
\end{lemma}
\begin{proof}
We first notice that $\C'\subseteq \D$ otherwise $\C'=\C$ implying $\C$ being symmetric. Then $\D''=\D$ by \cite[Theorem 3.10]{drinfeld2010braided}.

On  one hand, $\D=\C'\Rightarrow \D'=\C''=\C\vee\C'=\C$, where the second equality follows from \cite[Corollary 3.11]{drinfeld2010braided}. On the other hand, $\D'=\C\Rightarrow \D=\D''=\C'$.
\end{proof}

\begin{lemma}\label{integral1}
Let $\C=\C(\D,\kappa)$ be a braided MR fusion category. Assume that $\D\neq\C'$. Then $\C$ is weakly integral. In particular, if $\kappa\neq 0$ then $\C$ is integral.
\end{lemma}
\begin{proof}
We may assume that $\C$ is not symmetric. We notice that $\C'$ is a fusion subcategory of $\D'$ and $\D'\subseteq \D$ by Lemma \ref{centralizer}.

By \cite[Theorem 3.4]{drinfeld2010braided}, we have
$$\FPdim(\D)\FPdim(\D')=\FPdim(\C)\FPdim(\D\cap\C')=\FPdim(\C)\FPdim(\C'),$$
where the last equality follows from the fact that $\C'\subseteq\D$. Hence
\begin{equation}\label{equ-1}
\begin{split}
\frac{\FPdim(\C)}{\FPdim(\D)}=\frac{\FPdim(\D')}{\FPdim(\C')}.
\end{split}
\end{equation}

The right hand side is an integer because $\FPdim(\D')$ and $\FPdim(\C')$ are both integers, and $\FPdim(\D')$ divides $\FPdim(\C')$. Hence the left hand side is  also an integer. Let $a=\FPdim(\D)$. Then $\FPdim(\C)=2a+\frac{\kappa^2+\kappa\sqrt{\kappa^2+4a}}{2}$ is an integer and hence $\C$ is weakly integral. If $\kappa\neq0$ then $\sqrt{\kappa^2+4a}$ is an integer. It follows that $\FPdim(X_n)=\frac{\kappa+\sqrt{\kappa^2+4a}}{2}$ is an integer since it is a rational algebraic integer. Hence $\C$ is integral in this case.
\end{proof}

\begin{remark}
The assumption $\D\neq\C'$ is necessary. In fact, if  $\D=\C'$ then  the right side of  Equation (\ref{equ-1}) is $\frac{\FPdim(\C)}{\FPdim(\D)}$ and therefore  Equation (\ref{equ-1}) produces nothing. One counterexample is a Fibonacci category which is of the form $\C(\vect,1)$. It is non-degenerate and has Frobenius-Perron dimension $5+\sqrt{5}$.
\end{remark}

\begin{lemma}\label{integral2}
Let $\C=\C(\D,\kappa)$ be a braided MR fusion category. Assume that $\D=\C'$. Then $\C$ is one of the following:

(1)\, $\C$ is an equivariantization of a pointed fusion category.

(2)\, $\C$ is an equivariantization of a Fibonacci category $\F$.

(3)\, $\C$ is an equivariantization of the category $\svect\boxtimes\F$.
\end{lemma}
\begin{proof}
Let $\D=\C'=\Rep(A,u)$ for some finite group $A$ and a central element $u\in A$. Assume first that $u\ne 1$.

\medbreak
Let $\Rep(H):=\Rep(A/\langle u\rangle )\subseteq \C'$ be the  maximal Tannakian subcategory of $\Rep(A,u)$. By \cite[Theorem 4.18(ii)]{drinfeld2010braided}, there exists a braided fusion category $\B$ and an action of $H$ on $\B$ such that $\B^H=\C$ and $\vect^H=\Rep(H)$. Moreover, \cite[Proposition 4.56]{drinfeld2010braided} shows that $\B$ is slightly degenerate.On the other hand, $\Rep(H)$ is also a fusion subcategory of $\C'$, hence there exists a category $\B_1\subset \B$ such that $\B_1^H=\C'$. The Frobenius-Perron dimension of $\B_1$ is $\frac{\FPdim(\C')}{\FPdim(H)}=2$ by \cite[Proposition 4.26]{drinfeld2010braided}. Notice that $\B_1=\svect$ is the M\"{u}ger center of $\B$.

Let $\mathcal{O}_1,\mathcal{O}_2,\cdots,\mathcal{O}_s$ be the orbits of the simple objects of $\B$ under the action of $H$. Without loss of generality, we may assume that $\mathcal{O}_1=\{1\},\mathcal{O}_2=\{\delta\}$, where $\Irr(\B_1)=\{1,\delta\}$. Since $\B_1^H=\C'$ and there is only one simple object of $\C$ not contained in $\C'$, there are exactly three orbits. Assume that $\mathcal{O}_3=\{ Y_1,Y_2,\cdots, Y_t\}$. Let $d=\FPdim(Y_i)$. Then $Y_1\otimes Y_1^*=\1\oplus\sum_{i=1}^{t}a_iY_i$. Notice that $\delta$ can not appear in the decomposition of $Y_1\otimes Y_1^*$ by \cite[Proposition2.6]{etingof2011weakly}. Counting Frobenius-Perron dimensions on both sides, we get $d^2=1+d\sum_{i=1}^{t}a_i$ which shows that $d=1$ or $d=\frac{1+\sqrt{5}}{2}$. If $d=1$ then $\B$ is pointed.

If $d=\frac{1+\sqrt{5}}{2}$ then we will show that $t=2$. In fact, if $t=1$ then $\delta\otimes Y_1=Y_1$ which contradicts \cite[Proposition2.6]{etingof2011weakly}. If $t\geqslant 3$ then we may assume $\delta \otimes Y_1=Y_2$ and hence $\delta \otimes Y_2=Y_1$. It follows that  simple object $\delta\otimes Y_3$ can not be isomorphic to $Y_1$ which contradicts the fact that $H$ acts transitively on $\mathcal{O}_3$. In this case, $\B\cong \svect\boxtimes \F$, where $F$ is a Fibonacci category.

\medbreak
Assume now that $u=1$. Then $\Rep(A)=\C'$ is Tannakian subcategory of $\C$. Again by \cite[Theorem 4.18(ii)]{drinfeld2010braided}, there exists a braided fusion category $\B$ and an action of $H$ on $\B$ such that $\B^H=\C$ and $\vect^H=\Rep(A)$. Let $\mathcal{O}_1,\mathcal{O}_2,\cdots,\mathcal{O}_s$ be the orbits of the simple objects of $\B$ under the action of $H$, where $\mathcal{O}_1=\{1\}$. Then $s=2$ since $\vect^H=\Rep(A)$ and $\Irr(\Rep(A))\cup \{X_n\}=\Irr(\C)$. Assume that $\mathcal{O}_2=\{ Y_1,Y_2,\cdots, Y_t\}$.

If $t\geqslant 2$ then there exists $1\leqslant i\leqslant t$ such that $Y_i^*\neq Y_1$. Let $d=\dim(Y_i)$. Then $Y_1\otimes Y_i=\sum_{i=1}^{t}a_iY_i$. Counting categorical dimensions on both sides, we get $d^2=d\sum_{i=1}^{t}a_i$ which shows that $d$ is an integer. By the proof of \cite[Proposition 8.22] {etingof2005fusion}, $\frac{\dim(\B)}{\FPdim(\B)}=d$. This shows that $d$ divides $\dim(\B)$. On the other hand, $\dim(\B)=1+td$. So we have $d=1$. This proves that $\B$ is pointed.

If $t=1$ then $\B$ is a fusion category of rank $2$. By the classification of rank $2$ fusion categories, $\B$ is either pointed, or a Fibonacci category. This completes the proof.
\end{proof}

\begin{remark}
In \cite{etingof2011weakly}, the authors introduced the notion of a weakly  group-theoretical fusion category. The class of weakly group-theoretical categories is closed under taking equivariantizations \cite[Proposition 4.1]{etingof2011weakly}. Hence a braided MR fusion category with $\D=\C'$ is weakly  group-theoretical. In particular, if $\C$ fits into (1) of Lemma \ref{integral2} then it is group-theoretical by \cite[Theorem 7.2]{naidu2009fusion}.
\end{remark}

\medbreak

Combining Lemma \ref{integral1} and Lemma \ref{integral2}, we get the following result.

\begin{theorem}\label{integral3}
Let $\C=\C(\D,\kappa)$ be a braided MR fusion category. Then $\C$ is one of the following:

(1)\, $\C$ is weakly integral;

(2)\, $\C$ is an equivariantization of a Fibonacci category $\F$;

(3)\, $\C$ is an equivariantization of the category $\svect\boxtimes\F$.
\end{theorem}

\subsection{Symmetric MR fusion categories}
To classify symmetric fusion categories, it suffices to study the categories of representations of finite groups, see \cite{deligne1990categories}.

\medbreak
Let $\C=\Rep(G,u)$ be a  symmetric category and $\D$ be a  symmetric subcategory of $\C$. Assume that $\Irr(\C) = \{1 = \chi_1,\chi_2, \cdots ,\chi_k \}$, $\Irr(\D) =\{\chi_1 ,\chi_2 , \cdots ,\chi_{k - 1} \}$.

\medbreak

In \cite{1983Characters}, Gagola characterized finite groups $G$ which have an irreducible character $\chi$ such
that $\chi$ does not vanish on exactly two conjugacy classes of $G$.  Now assume that $|G|>2$.  Our next result shows that our context is equivalent to the context of Gagola.

\begin{theorem}\label{finitegroup}
	Keep the notation above. Then the following two statements are equivalent:
	
	(1)\, ${\Irr}(\C)$ has a maximal symmetric subcategory $\D$;
	
	(2)\, There exists $\chi\in{\Irr}(\C)$ such that   $\chi$ does not vanish on exactly two conjugacy classes of $G$.	
	
\end{theorem}

\begin{proof}
	First,  assume statement $(1)$ holds. As $\D$ is a subcategory, it is well-known that there exists a normal subgroup $N$ of $G$ such that  $\D=\Irr(G/N) =\{\chi_1 ,\chi_2 , \cdots ,\chi_{k - 1} \}$ for a normal subgroup $N$ of $G$ (see \cite[p139]{blau1991}). Obviously, $N$ is the intersection of the kernels of $\chi_j$ with $1\le j\le n-1$ and   $\Irr(G) =\Irr(G/N)\cup\{\chi_k\}$.	
	Note that
	$$
	|{\rm Cla}(G/N)|=|\Irr(G/N)|=|{\rm Cla}(G)|-1,
	$$
	where ${\rm Cla}(G/N)$ and ${\rm Cla}(G)$ are respectively the set of conjugacy classes of $G/N$ and $G$.
	Hence, one can easily see that nonidentity elements of $N$ consists of a single conjugacy class $\C_2$ of $G$.
	
	Since $N\subseteq {\rm ker}(\chi_i)$, $1\leqslant i\leqslant k-1$, the character table of $G$ looks likes the following:
	\begin{center}
		\begin{tabular}{c|ccccccc}
			$\cdot$ & $\{1\}=\C_1$ & $\C_2$ & $\C_3$ &$\cdots$ & $\C_i$ & $\cdots$ & $\C_k$ \\
			\hline
			$\chi_1$ & $n_1=1$ & $n_1$ & $c_{13}$ &$\cdots$ & $c_{1i}$ & $\cdots$ & $c_{1k}$\\
			$\chi_2$ & $n_2$ & $n_2$ & $c_{23}$ &$\cdots$ & $c_{2i}$ & $\cdots$ & $c_{2k}$\\
			$\vdots$ & $\vdots$ & $\vdots$ &$\vdots$ & & $\vdots$ & $ $ & $\vdots$ \\
			$\chi_{k-1}$ & $n_{k-1}$ & $n_{k-1}$ &$c_{k-1,3}$ & $\cdots$ & $c_{k-1,i}$ & $\cdots$ & $c_{k-1,k}$\\
			$\chi_{k}$ & $n_{k}$ & $c_{k2}$ &$c_{k3}$ & $\cdots$ & $c_{k,i}$ & $\cdots$ & $c_{k,k}$
		\end{tabular}
	\end{center}
	
	Obviously, the complex conjugate $\overline{c_{k2}}$ is not equal to $ n_k$ as the character table is nonsingular. For any $i\ge 2$, if we choose $x_i\in\C_i$, then by Second Orthogonality Relation of irreducible characters of finite group (see \cite[Theorem 2.18]{charactertheory}),  for any $i>2$ we have
	\begin{equation}
	\begin{split}
	&\sum_{j=1}^{k}\overline{\chi_j(x_2)}\chi_j(x_i)=0=n_1\chi_1(x_i)+\cdots+n_{k-1}\chi_{k-1}(x_i)+\overline{c_{k2}}\chi_k(x_i),\\
	&\sum_{j=1}^{k}\chi_j(1)\chi_j(x_i)=0=n_1\chi_1(x_i)+\cdots+n_{k-1}\chi_{k-1}(x_i)+n_{k}\chi_k(x_i).
	\end{split}\nonumber
	\end{equation}
	
	Hence $(\overline{c_{k2}}-n_k)\chi_k(x_i)=0$. As $\overline{c_{k2}}-n_k\neq0$, we have $\chi_k(x_i)=0$, $i=2,\ldots, k$. Moreover, $c_{k2}\ne 0$. Thus, statement $(2)$ holds.
	\medbreak
	
	Conversely, assume statement $(2)$ holds. By \cite[Lemma 2.1]{1983Characters}, the irreducible character $\chi$ in statement $(2)$ is the unique faithful irreducible character of $G$. Moreover,   $G$ possesses a unique minimal normal subgroup $N$ such that $N$ is an elementary abelian $p$-group for a prime $p$ and the nonidentity element of $N$  consists of a conjugacy class of $G$. It then follows that for any $\varphi\in {\Irr}(\C)\setminus \{\chi\}$, $N\subseteq {\rm ker}(\varphi)$. This yields that
	$$
	\D:={\Irr}(\C)\setminus \{\chi\}={\Irr}(G/N)
	$$
	is a subcategory of $\C$ with maximal rank. The proof of the theorem is complete.
\end{proof}

\subsection{Non-degenerate MR fusion categories}
\begin{lemma}\label{prop1}
Let $\C(\D,\kappa)$ be an MR fusion category (not necessarily braided). Assume $\C$ is weakly integral and $\FPdim(X_n)^2$ divides $\FPdim(\C)$ then  $\D=\C_{ad}$ and $\U(\C)=\mathbb{Z}_2$.
\end{lemma}
\begin{proof}
We may assume that there exists $m$ such that $m(\FPdim(X_n)^2)=\sum_{i=1}^n\FPdim(X_i)^2=\FPdim(\C)$. On the other hand, the fusion rules of $\C$ show that $\FPdim(X_n)^2=\sum_{i=1}^{n-1}\FPdim(X_i)^2+\kappa\FPdim(X_n)$. Then $(m-2)\FPdim(X_n)^2=-\kappa\FPdim(X_n)$, which shows that $\kappa=0$. The fusion rules of $\C$ then imply that $\C_{ad}=\D$ and hence $\U(\C)=\mathbb{Z}_2$.
\end{proof}

\begin{corollary}
Let $\C(\D,\kappa)$ be an MR fusion category. Assume that $\C$ is braided and nilpotent. Then $\FPdim(\C)$ is a power of $2$.
\end{corollary}

\begin{proof}
Since $\C$ is nilpotent, it is weakly integral and hence we may assume that $\FPdim(\C)=p_1^{t_1}p_2^{t_2}\cdots p_s^{t_s}$, where $p_1,p_2,\cdots p_s$ are distinct primes and $t_i>0$. Then $\C$ has a unique decomposition of Delige product $\C=\C_{p_1}\boxtimes \C_{p_2}\boxtimes\cdots \boxtimes\C_{p_s}$ \cite[Theorem 1.1]{drinfeld2007g}, where each $\C_{p_i}$ is a fusion subcategory of $\C$ of dimension $p_i^{t_i}$.
Without loss of generality, we may assume that $X_n$ belongs to $\C_{p_1}$. Then the decomposition of $X_n\otimes X_n^*$ shows that all simple objects of $\C$ belong to $\C_{p_1}$. Hence $\C=\C_{p_1}$ and $\FPdim(\C)$ is a power of $p_1$.

By \cite[Corollary 5.3]{gelaki2008nilpotent}, $\FPdim(X_n)^2$ divides $\FPdim(\C_{ad})$ and hence divides $\FPdim(\C)$. It follows from Lemma \ref{prop1} that $\FPdim(\C)$ is even. Thus $p_1=2$ and $\FPdim(\C)$ is a power of $2$.
\end{proof}

Recall that a braided fusion category is called prime if it contains no proper non-trivial non-degenerate fusion subcategories.
\begin{proposition}\label{prime}
Let $\C=\C(\D,\kappa)$ be a braided MR fusion category. Then $\C$ is prime.
\end{proposition}
\begin{proof}
Assume on the contrary that $\C$ contains a proper non-degenerate fusion subcategory $\A$.Then $\A$ must be a fusion subcategory of $\D$. By M\"{u}ger  Decomposition Theorem \cite[Theorem 3.13]{drinfeld2010braided}, we have tensor equivalences $\C=\A\boxtimes \B$ and $\D=\A\boxtimes \E$,  where $\B$ and $\E$ are respectively the centralizer of $\A$ in $\C$ and $\D$. Comparing their ranks, we have $n=\rk(\A)\rk(\B)$ and $n-1=\rk(\A)\rk(\E)$ which deduces that $1=\rk(\A)(\rk(\B)-\rk(\E))$. This is impossible.
\end{proof}

A non-degenerate braided fusion category $\C$ is called a minimal extension of a braided fusion
category $\A$ if $A\subset \C$ and $\mathcal{Z}_2(\A)=\A'$.
\medbreak

\begin{theorem}\label{non-degenerate}
Let $\C=\C(\D,\kappa)$ be a non-degenerate MR fusion category. Then $\C$ is exactly one of the following:

(1)\, a pointed modular category $\C(\mathbb{Z}_2,\pm i)$;

(2)\, a Fibonacci category;

(3)\, a minimal extension of a slightly degenerate fusion category.
\end{theorem}
\begin{proof}
If $\D=\C'=\vect$ then $\C$ is a rank two fusion category. Then (1) and (2) follow from the classification of rank $2$ fusion categories \cite{ostrik2003fusion} and pointed modular categories \cite[Example 5.1]{drinfeld2007g}.

\medbreak
In the other case, $\C'\subseteq\D$ is a fusion subcategory of $\D$. Hence $\C$ is weakly integral by Lemma \ref{integral1}.

Let $X$ be a simple object of $\C$. Then $\FPdim(X)^2$ divides $\FPdim(\C)$ by \cite[Theorem 2.11]{etingof2011weakly}. Then $\D=\C_{ad}$ and $\U(\C)=\mathbb{Z}_2$ by Lemma \ref{prop1}. Hence $\C_{pt}=\vect_{\mathbb{Z}_2}^{\omega}$ for some $3$-cocycle $\omega\in H^3(\mathbb{Z}_2,k^{\times})$ by \cite[Theorem 6.2]{gelaki2008nilpotent}.

It follows from \cite[Corollary 3.27]{drinfeld2010braided} that $\D'=\C_{ad}'=\C_{pt}$. On the other hand, $\C_{pt}$ is contained in $\D$. Hence $\Z_2(\D)=\D\cap \D'=\C_{pt}$. Then $\C_{pt}$ is a symmetric fusion subcategory and thus it is either a Tannakian subcategory, or a category $\svect$ of super vector spaces. We shall prove that $\C_{pt}$ is not Tannakian and hence $\D$ is a slightly degenerate fusion subcategory.

Suppose on the contrary that $\C_{pt}$ is Tannakian. Set $G(\C_{pt})=\{\1,g\}$, $\Irr(\C) = \{1 = X_1 ,X_2 , \cdots ,X_n \}$, $\Irr(\D) =
\{X_1 ,X_2 , \cdots ,X_{n - 1} \}$. Then $g\otimes X_n=X_n$ by the fusion rules of $\C$. The entry $s_{g,X_n}$ of the $\mathcal{S}$-matrix is
\begin{equation}
\begin{split}
s_{g,X_n}&=\theta_{g}^{-1}\theta_{X_n}^{-1}\sum_{Z\in \Irr(\C)}N_{gX_n}^{Z}\theta_Zd_Z\\
&=\theta_{g}^{-1}\theta_{X_n}^{-1}\theta_{X_n}d_{X_n}\\
&=\theta_{g}^{-1}d_{X_n}\\
&=d_{X_n},
\end{split}
\end{equation}
where the last equality follows from the fact $\theta_{g}=1$ since $\C_{pt}$ is Tannakian.

Notice that, for all $i=1,\cdots, n-1$, $X_i$ is contained in $\C_{ad}$ and $g$ is contained in $\C_{ad}'$. Hence $s_{g,X_i}=d_{X_i}$ by \cite[Proposition2.5]{muger2003structure}. Thus we get that the rows of the $\mathcal{S}$-matrix corresponding to $g$ and $1$ are equal. This implies that the $\mathcal{S}$-matrix is degenerate which contradicts the non-degeneracy  of $\C$. This completes the proof.
\end{proof}

\begin{corollary}
Let $\C=\C(\D,\kappa)$ be a non-degenerate MR fusion category. If $\D$ is pointed then $\C$ is one of the following:

(1)\, a pointed modular category $\C(\mathbb{Z}_2,\pm i)$;

(2)\, a Fibonacci category;

(3)\, an Ising category.
\end{corollary}
\begin{proof}
It suffices to consider the third case. If $\D$ is bigger than $\svect$ then there exists a proper non-degenerate fusion subcategory $\A$ such that $\D=\A\boxtimes \svect$. This contradicts Proposition \ref{prime}. Hence $\D=\svect$ and $\C$ is an Ising category.
\end{proof}

\section{MR fusion categories with rank $\leq 4$}\label{sec6}
Let $\C(\D,\kappa)$ be an MR fusion category. As we have seen, it is impossible to
classify $\C(\D,\kappa)$ completely, because it at least contains the classification
of category of representations of a semisimple quasi-Hopf algebra. In this section, we aim to classify $\C(\D,\kappa)$  with rank $\leq 4$.

\medbreak
According to the conclusion of \cite{ostrik2003fusion}, there are only two classes of fusion categories with rank $2$: pointed fusion categories and Fibonacci categories. The pointed fusion category with only two simple objects is obviously a category of the form of $\C(\vect,0)$, where $\vect $ is the trivial fusion category. The fusion rules of a Fibonacci category is $X = \1 \oplus X $, where $\1$ is the unit object and $X $ is the unique noninvertible object. A Fibonacci category is of the form $\C(\vect,1)$.

\medbreak
If $\C(\D,\kappa)$ has rank $3$, then $\D$ is an integral fusion category with rank $2$.
From the above discussion, we know that $\D$ is a pointed fusion category, so $\C$ is a near-group category. More precisely, $\C$ is an Ising category, or the category of representations of the group $S_3$ or its  twisted version, see \cite[Theorem 1.1]{2013ostrikpivotal}.

\medbreak
In the rest of this section, we are devoted to the classification of fusion category $\C(\D,\kappa)$ with rank $4$.

\begin{lemma}\label{lemmaS_3}
The fusion ring below can not  be categorizable, where $\kappa$ is non-negative.
\[
X X = 1 +X + Y,
\quad
YY = 1,
\quad
X Y = Y X = X,
\]
\[
X Z = Z X = 2Z,
YZ = Z Y = Z,
Z Z = 1 + 2X +Y +\kappa Z.
\]
\end{lemma}
\begin{proof}
Assume on the contrary there exists a fusion category $\C$ which has such fusion rules. Then $\{1,X,Y\}$ generates a fusion subcategory $\D$ of $\C$.  Moreover, $\FPdim(X)=2$ and $\FPdim(Y)=1$. Applying the homomorphism $\FPdim$  on both sides of the last formula, we can get that $\FPdim (Z)=\frac{\kappa + \sqrt {\kappa^2 + 24} }{2}$.

Considering the left tensor product of $X,Y $ and $Z $ acting on the basis $\Irr (\C)=\{\1, X, Y, Z \}$ of $K (\C) $,  we obtain three matrices

\[
M_X = \left( {{\begin{array}{*{20}c}
0 \hfill & 1 \hfill & 0 \hfill & 0 \hfill \\
1 \hfill & 1 \hfill & 1 \hfill & 0 \hfill \\
0 \hfill & 1 \hfill & 0 \hfill & 0 \hfill \\
0 \hfill & 0 \hfill & 0 \hfill & 2 \hfill \\
\end{array} }} \right),
\quad
M_Y = \left( {{\begin{array}{*{20}c}
0 \hfill & 0 \hfill & 1 \hfill & 0 \hfill \\
0 \hfill & 1 \hfill & 0 \hfill & 0 \hfill \\
1 \hfill & 0 \hfill & 0 \hfill & 0 \hfill \\
0 \hfill & 0 \hfill & 0 \hfill & 1 \hfill \\
\end{array} }} \right),
\quad
M_Z = \left( {{\begin{array}{*{20}c}
0 \hfill & 0 \hfill & 0 \hfill & 1 \hfill \\
0 \hfill & 0 \hfill & 0 \hfill & 2 \hfill \\
0 \hfill & 0 \hfill & 0 \hfill & 1 \hfill \\
1 \hfill & 2 \hfill & 1 \hfill & \kappa \hfill \\
\end{array} }} \right).
\]

Set $M =I_4 + M_X^2 + M_Y^2 + M_Z ^ 2 $, where $I_4 $ is the $4\times
4$ identity matrix. The four eigenvalues of $M$  are $f_1 = \frac{1}{2}(24 + \kappa^2 + \kappa\sqrt {24 +
\kappa^2} )$, $f_2 = \frac{1}{2}(24 + \kappa^2 - \kappa\sqrt {24 + \kappa^2} )$, $f_3 = 3$ and $f_4 =2$. These eigenvalues are called codegrees of $\C$  in \cite{2013ostrikpivotal}.

Since $K (\C) $ is a commutative ring, it has four $1$-dimensional representations.
The object $\I (\1) $ in Drinfeld center $\Z (\C) $ is a direct sum of four simple objects, and the multiplicity of each simple object is $1$, see  \cite[Theorem 2.13]{2013ostrikpivotal}.
Since $\I (\1) $ is a commutative algebra in $\Z (\C) $, the unit object $1$ is a simple object in $\I (\1)$.
Therefore, we can assume $\I (\1) = \1 \oplus A \oplus B\oplus E$.

From \cite[Theorem 2.13]{2013ostrikpivotal}, we get
\begin{equation}
\begin{split}
\FPdim (A) &= \frac {f_1 }{f_2 } =
\frac{1}{12}(12 + \kappa^2 + \kappa\sqrt {24 + v^2} ),\\
\FPdim(B) &= \frac{f_1}{f_3 } = \frac{1}{6}(24 + \kappa^2 + \kappa\sqrt {24 + \kappa^2} ),\\
\FPdim(E) &=\frac{f_1 }{f_4 } = \frac{1}{4}(24 + \kappa^2 + \kappa\sqrt {24 + \kappa^2}).
\end{split}
\end{equation}

From \cite[Proposition 5.4]{etingof2005fusion},
\[
\F(\I(\1)) = F(\1) \oplus \F(A) \oplus \F(B) \oplus \F(E) = \mathop \oplus
\limits_{T \in \mbox{Irr}(C)} T \otimes T^\ast = 4 \cdot \1 \oplus 3X \oplus
2Y \oplus \kappa Z.
\]

We assume
\begin{equation}\label{decom01}
\begin{split}
\F(A) = a_1 X + a_2 Y + a_3 Z,\\
\F(B) = b_1 X + b_2 Y + b_3 Z,\\
\F(E) = c_1 X + c_2 Y + c_3 Z.
\end{split}
\end{equation}

Then
\begin{equation}\label{coeff01}
\begin{split}
a_1 +b_1+c_1=3, a_2 +b_2+c_2=2, a_3 +b_3+c_3=\kappa.
\end{split}
\end{equation}

Considering Frobenius-Perron dimensions on both sides of (\ref{decom01}), we get

\begin{gather}
(\kappa - 6a_3 )\sqrt {24 + \kappa^2} = 24a_1 + 12a_2 + 6a_3 \kappa - \kappa^2,\label{eqr2}\\
(\kappa - 3b_3 )\sqrt {24 + \kappa^2} = 2b_1 + 6b_2 + 3b_3 \kappa -18- \kappa^2,\label{eqr3}\\
(\kappa - 2c_3 )\sqrt {24 + n^2} = 8c_1 + 4c_2 + 2c_3 \kappa -20- \kappa^2.\label{eqr4}
\end{gather}

We first prove that $\sqrt {24 + \kappa^2}$ is an integer. Assume on the contrary that $\sqrt {24 + \kappa^2}$ is an irrational number.  Then equations (\ref{eqr2}-\ref{eqr4}) induce that
\begin{equation}\label{coeff02}
\begin{split}
&\kappa - 6a_3=0, \kappa - 3b_3=0, \kappa - 2c_3=0,\\
&2a_1 + a_2=0, b_1 + 3b_2  -9=0,2c_1 + c_2 -5=0.
\end{split}
\end{equation}

Together with equations in(\ref{coeff01}), it is easy to check that the above system of equations has no solutions. This proves that $\sqrt {24 + \kappa^2}$ is an integer, which shows that $\kappa = 1$ or $5$.

If $\kappa=1$, then equation (\ref{eqr2}) is

\[
5(1 - 6a_3 ) = 24a_1 + 12a_2 + 6a_3 - 1.
\]

If $a_1 = a_2 = a_3 = 0 $, then the left side of the above equation is $5$ and the right side is $-1$,
a contradiction. If $a_3\neq 0$ then the left side is negative and the right side is positive, also a contradiction. If $a_3=0$, $a_1$ or $a_2$ is not $0$ then the left side is $5$ and the right side is great or equal to $11$, still a contradiction.

If $\kappa = 5 $, then equation (\ref{eqr2}) is $5 = - 2 a_1 + a_2 + 6a_3 $. Combining with equations in (\ref{coeff01}), we get $a_1 = 2, a_2 = 1, a_3 = 0$. Taking this result into equations (\ref{coeff01}), we get
\begin{center}
$b_1 = 0, b_2 = 1, b_3 = 0$ and $c_1 = 1, c_2 = 0, a_3 = 3$.
\end{center}

Again by \cite[Proposition 5.4]{etingof2005fusion}, we have
\begin{equation}
\begin{split}
\F(\I(X)) = 3 \cdot \1 \oplus 9X \oplus 3Y \oplus 10Z,\\
\F(\I(Y)) = 2 \cdot \1 \oplus 3X \oplus 4Y \oplus 5Z,\\
\F(\I(X)) = 5 \cdot \1 \oplus 10X \oplus 5Y \oplus 37Z.
\end{split}
\end{equation}

From
\begin{equation*}
\begin{split}
\dim\Hom(\I(1),\I(Y))&=\dim\Hom(\F(\I(1)),Y)=2,\\
\dim\Hom(A,\I(Y))&=\dim\Hom(\F(A),Y)=1,\\
\dim\Hom(B,\I(Y))&=\dim\Hom(\F(B),Y)=1,\\
\dim\Hom(E,\I(Y))&=\dim\Hom(\F(E),Y)=0,\\
\dim\Hom(\I(Y),\I(Y))&=\dim\Hom(\F(\I(Y)),Y)=4,
\end{split}
\end{equation*}
we can write $\I(Y)=A\oplus B\oplus G_1\oplus G_2$, where $G_1,G_2$ are non-isomorphic simple objects of $\Z(\C)$ which are both different from $E$. So we have $$\F(\I(Y))=\F(A)\oplus \F(B)\oplus \F(G_1)\oplus \F(G_2)=2\cdot\1\oplus 3X\oplus 4Y\oplus 5Z.$$
Hence $\F(G_1)\oplus \F(G_2)=X\oplus 2Y\oplus 3Z$. In particular, $Z$ must appear in $\F(G_1)$ or $\F(G_2)$ as a summand.

From
\begin{equation*}
\begin{split}
\dim\Hom(\I(1),\I(Z))&=\dim\Hom(\F(\I(1)),Z)=5,\\
\dim\Hom(A,\I(Z))&=\dim\Hom(\F(A),Z)=0,\\
\dim\Hom(B,\I(Z))&=\dim\Hom(\F(B),Z)=2,\\
\dim\Hom(E,\I(Z))&=\dim\Hom(\F(E),Z)=3,\\
\dim\Hom(\I(Z),\I(Z))&=\dim\Hom(\F(\I(Z)),Z)=37,
\end{split}
\end{equation*}
we can write $\I(Z)=2B\oplus 3E\oplus Q$, where $Q$ does not contain $A$ as a summand.

It is easy to see that $\dim\Hom(\I(Y),I(Z))=\dim\Hom(\F(\I(Y)),Z)=5$. But the left hand side is
$\dim\Hom(A\oplus B\oplus G_1\oplus G_2,2B\oplus 3E\oplus Q)=5+\dim\Hom(G_1\oplus G_2,Q)$ which shows that $\dim\Hom(G_1\oplus G_2,Q)=0$. This shows that $Q$ and hence $\I(Z)$ can not contain $G_1$ or $G_2$ as a summand. It follows that $\dim\Hom(G_i,\I(Z))=\dim\Hom(\F(G_i),Z)=0$ which contradicts the conclusion obtained above. This finishes the proof.
\end{proof}

Since the proof of the Lemma below follows the same line of Lemma \ref{lemmaS_3}, we put it in the Appendix.
\begin{lemma}\label{labelZ_3}
If a fusion category $\C$ has the following fusion ring  then $\kappa=2$ or $\kappa$ is divisible by $3$.
\[
X  X = Y,
\quad
Y  Y = X,
\quad
X Y = Y X = \1,
\]
\[
X  Z = Z  X = Z,
Y Z = Z Y = Z,
Z  Z = 1 + X + Y + \kappa Z.
\]
\end{lemma}

\begin{theorem}\label{rank4}
Let $\C=\C(\D,\kappa)$ be an MR fusion category. Assume that the rank of $\C$ is $4$. Then $\C$ is a near-group fusion category and $\kappa=2$ or $\kappa$ is divisible by $3$.
\end{theorem}

\begin{proof}
Under our assumption, $\D$ is an integral fusion category of rank $3$ by Proposition \ref{Prop3.2}. By the classification of rank $3$ pivotal fusion categories, $\D$ is either pointed, or the category of representations of the group $S_3$ or its  twisted version, see \cite[Theorem 1.1]{2013ostrikpivotal}. If $\D$ is  the category of representations of the group $S_3$ or its  twisted version, then $\C$ has the following fusion rules:

\[
X X = 1 +X + Y,
\quad
YY = 1,
\quad
X Y = Y X = X,
\]
\[
X Z = Z X = 2Z,
YZ = Z Y = Z,
Z Z = 1 + 2X +Y +\kappa Z.
\]
But Lemma \ref{lemmaS_3} shows that it is impossible.

\medbreak
If $\D$ is pointed then $\C$ has the following fusion rules:
\[
X  X = Y,
\quad
Y  Y = X,
\quad
X Y = Y X = \1,
\]
\[
X  Z = Z  X = Z,
Y Z = Z Y = Z,
Z  Z = 1 + X + Y + \kappa Z.
\]

The theorem then follows from Lemma \ref{labelZ_3}.
\end{proof}

\section{Appendix}
\begin{proof}[Proof of Lemma \ref{labelZ_3}]
It is obvious that  $\{1,X,Y\}$ generates a pointed fusion subcategory $\D$ of $\C$.

Considering the left tensor product of $X,Y $ and $Z $ acting on the basis $\Irr (\C)=\{\1, X, Y, Z \}$ of $K (\C) $,  we obtain three matrices

\[
M_X = \left( {{\begin{array}{*{20}c}
0 \hfill & 0 \hfill & 1 \hfill & 0 \hfill \\
1 \hfill & 0 \hfill & 0 \hfill & 0 \hfill \\
0 \hfill & 1 \hfill & 0 \hfill & 0 \hfill \\
0 \hfill & 0 \hfill & 0 \hfill & 1 \hfill \\
\end{array} }} \right),
\quad
M_Y = \left( {{\begin{array}{*{20}c}
0 \hfill & 1 \hfill & 0 \hfill & 0 \hfill \\
0 \hfill & 0 \hfill & 1 \hfill & 0 \hfill \\
1 \hfill & 0 \hfill & 0 \hfill & 0 \hfill \\
0 \hfill & 0 \hfill & 0 \hfill & 1 \hfill \\
\end{array} }} \right),
\quad
M_Z = \left( {{\begin{array}{*{20}c}
0 \hfill & 0 \hfill & 0 \hfill & 1 \hfill \\
0 \hfill & 0 \hfill & 0 \hfill & 1 \hfill \\
0 \hfill & 0 \hfill & 0 \hfill & 1 \hfill \\
1 \hfill & 1 \hfill & 1 \hfill & \kappa \hfill \\
\end{array} }} \right).
\]

Set $M =I_4 + M_X^2 + M_Y^2 + M_Z ^ 2 $, where $I_4 $ is the $4\times
4$ identity matrix. The four eigenvalues of $M$  are $f_1 = 3$, $f_2 =3$, $f_3 = \frac{1}{2}(12 + \kappa ^2 - \kappa\sqrt {12 +
\kappa ^2} )$ and $f_4 = \frac{1}{2}(12 + \kappa ^2 + \kappa \sqrt {12 + \kappa ^2} )$.

If $\sqrt {12 + \kappa ^2} $ is an integer then $\kappa =2$. We are done in this case. In the rest of the proof, we assume that $\sqrt {12 + \kappa ^2} $ is irrational.

\medbreak
By the same arguments as in Lemma \ref{lemmaS_3}, we can assume that $\I (\1) = \1 \oplus A \oplus B\oplus E$ and 
\begin{equation}
\begin{split}
\FPdim (A) &=
\frac{1}{6}(12 + \kappa^2 + \kappa \sqrt {12 + \kappa^2} ),\\
\FPdim(B) &=  \frac{1}{6}(12 + \kappa^2 + \kappa\sqrt {12 + \kappa^2} ),\\
\FPdim(E) &=1+ \frac{1}{6}( \kappa^2 +  \kappa\sqrt {12 + \kappa^2}).
\end{split}
\end{equation}

From \cite[Proposition 5.4]{etingof2005fusion}, we have
\[
\F(\I(\1)) = \F(\1) \oplus \F(A) \oplus \F(B) \oplus \F(E) = \mathop \oplus
\limits_{T \in \mbox{Irr}(C)} T \otimes T^\ast = 4 \cdot \1 \oplus X \oplus
Y \oplus \kappa Z.
\]

We assume
\begin{equation}\label{decom001}
\begin{split}
\F(A) = a_1 X + a_2 Y + a_3 Z,\\
\F(B) = b_1 X + b_2 Y + b_3 Z,\\
\F(E) = c_1 X + c_2 Y + c_3 Z.
\end{split}
\end{equation}

Then
\begin{equation}\label{coeff001}
\begin{split}
a_1 +b_1+c_1=1, a_2 +b_2+c_2=1, a_3 +b_3+c_3=\kappa.
\end{split}
\end{equation}

Considering Frobenius-Perron dimensions on both sides of $\F(A)$, we get

\begin{gather}
(\kappa - 3a_3 )\sqrt {12 + \kappa^2} = -6a_1 + 6a_2 + 3a_3 \kappa - \kappa^2.
\end{gather}

Since $\sqrt {12 + \kappa^2}$ is irrational, $\kappa - 3a_3$ must be $0$. This implies that $\kappa$ is divisible by $3$.
\end{proof}

\section*{Acknowledgements}
The research of the first  author is partially supported by the Natural Science Foundation of Jiangsu Providence (Grant No. BK20201390) and the startup foundation for introducing talent of NUIST (Grant No. 2018R039). The research of the second author is partially supported by the Natural Science Foundation of China (Grant No. 11971189).


\begin{thebibliography}{10}

\bibitem{blau1991}
Z.~Arad,  H.~Blau.
\newblock On table algebras  applications to finite group theory.
\newblock {\em J. Algebra}, 138(1):137--185, 1991.

\bibitem{bruguieres2000categories}
A. Bruguieres.
\newblock Cat{\'e}gories pr{\'e}modulaires, modularisations et invariants des
  vari{\'e}t{\'e}s de dimension $3$.
\newblock {\em Math. Ann.}, 316(2):215--236, 2000.

\bibitem{deligne1990categories}
P. Deligne.
\newblock Cat{\'e}gories {T}annakiennes.
\newblock In {\em The Grothendieck Festschrift}, pages 111--195. Springer,
  1990.

\bibitem{dong2012frobenius}
J. Dong, S. Natale,  L. Vendramin.
\newblock {F}robenius property for fusion categories of small integral
  dimension.
\newblock {\em J. Algebra Appl.}, 14(2):1550011[17pages], 2015.

\bibitem{drinfeld2007g}
V. Drinfeld, S. Gelaki, D. Nikshych,  V. Ostrik.
\newblock Group-theoretical properties of nilpotent modular categories.
\newblock {\em preprint arXiv:0704.0195}, 2007.

\bibitem{drinfeld2010braided}
V. Drinfeld, S. Gelaki, D. Nikshych,  V. Ostrik.
\newblock On braided fusion categories {I}.
\newblock {\em Selecta Math., New Ser.}, 16(1):1--119, 2010.

\bibitem{egno2015}
P.~Etingof, S.~Gelaki, D.~Nikshych,  V.~Ostrik.
\newblock {\em Tensor Categories}.
\newblock Mathematical surveys  monographs, vol. 205. Amer. Math. Soc.,
  2015.

\bibitem{etingof2011weakly}
P. Etingof, D. Nikshych,  V. Ostrik.
\newblock Weakly group-theoretical  solvable fusion categories.
\newblock {\em Adv. Math.}, 226(1):176--205, 2011.

\bibitem{etingof2005fusion}
P. Etingof, D. Nikshych,  V. Ostrik.
\newblock On fusion categories.
\newblock {\em Ann. Math.}, 162(2):581--642, 2005.

\bibitem{1983Characters}
S.~M. Gagola.
\newblock Characters vanishing on all but two conjugacy classes.
\newblock {\em Pacific Journal of Mathematics}, 109(2):363--385, 1983.

\bibitem{gelaki2008nilpotent}
S. Gelaki, D. Nikshych.
\newblock Nilpotent fusion categories.
\newblock {\em Adv. Math.}, 217(3):1053--1071, 2008.

\bibitem{charactertheory}
I.~M. Isaacs.
\newblock {\em Character theory of finite groups}.
\newblock Academic Press,New York, 1976.

\bibitem{kassel1995quantum}
C. Kassel.
\newblock Quantum groups, {GTM} 155, 1995.

\bibitem{HONG20101000}
S. M. Hong, E. Rowell.
\newblock On the classification of the grothendieck rings of non-self-dual
  modular categories.
\newblock {\em Journal of Algebra}, 324(5):1000--1015, 2010.
\newblock Computational Algebra.

\bibitem{muger2003structure}
M. M{\"u}ger.
\newblock On the structure of modular categories.
\newblock {\em Proc. London Math. Soc.}, 87(02):291--308, 2003.

\bibitem{muger2004galois}
M. M{\"u}ger.
\newblock Galois extensions of braided tensor categories  braided crossed
  {G}-categories.
\newblock {\em J. Algebra}, 277(1):256--281, 2004.

\bibitem{naidu2009fusion}
D. Naidu, D. Nikshych,  S. Witherspoon.
\newblock Fusion subcategories of representation categories of twisted quantum
  doubles of finite groups.
\newblock {\em Internat. Math. Res. Notices}, 2009(22):4183--4219, 2009.

\bibitem{Nichols1996297}
W.~D. Nichols  M. B. Richmond.
\newblock The {G}rothendieck group of a {H}opf algebra.
\newblock {\em J. Pure Appl. Algebra}, 106(3):297 -- 306, 1996.

\bibitem{ostrik2003module}
V. Ostrik.
\newblock Module categories, weak {H}opf algebras  modular invariants.
\newblock {\em Transform. Groups}, 8(2):177--206, 2003.

\bibitem{2013ostrikpivotal}
V. Ostrik.
\newblock Pivotal fusion categories of rank 3.
\newblock {\em Mosc. Math. J.}, 15(2):373--396, 2015.

\bibitem{ostrik2003fusion}
V. Ostrik.
\newblock Fusion categories of rank 2.
\newblock {\em Math. Res. Lett.}, 10(2):177--183, 2003.

\bibitem{siehler2003near}
J. Siehler.
\newblock Near-group categories.
\newblock {\em Algebr. Geom. Topol.}, 3:719--775, 2003.

\bibitem{Tambara1998692}
D. Tambara, S. Yamagami.
\newblock Tensor categories with fusion rules of self-duality for finite
  abelian groups.
\newblock {\em J. Algebra}, 209(2):692 -- 707, 1998.

\end{thebibliography}

\end{document}